\def\IS{{\Bbb S}} 
 
\def\IK{{\Bbb K}}

\def\ID{{\Bbb D}}

\def\zbar{{\overline{z}}} 
 
\def\wbar{{\overline{w}}}

\documentclass[a4paper,  12pt,  times]{article}
\usepackage{secdot} 
\sectiondot{subsection} 
\usepackage{titling}
\usepackage{amsmath,  amsfonts,  amsthm, mathtools}
\newtheorem{theorem}{Theorem}
\newtheorem*{theorem*}{Theorem}
\newtheorem{corollary}{Corollary}
\newtheorem*{corollary*}{Corollary}

\newtheorem{lemma}{Lemma}
\newtheorem{conj}{Conjecture}

\title{On the Conformal Energy of Quasisymmetric and Quasim\"obius Mappings}
\author{Mihai Florincescu and Gaven Martin \thanks{Research of both authors supported by the NZ Marsden Fund, 
\newline New Zealand Institute for Advanced Study,  Massey University.
\newline
{\bf email:} m.florincescu@massey.ac.nz\newline
{\bf email:} g.j.martin@massey.ac.nz
\newline
\newline This article is based on the first author's MSc Thesis.
}}

\begin{document}
\maketitle
\begin{abstract}
This article identifies the conformal energy (or mean distortion) of extremal mappings of finite distortion with a given quasisymmetric mapping of the circle as boundary data. The conformal energy of $g_o:\IS\to\IS$ is
\begin{equation}\label{energy} \mathcal E(g_o)=-\frac{1}{4\pi^2}\iint_{\mathbb{S}\times \mathbb{S}}\log |g_o(\zeta)-g_o(\eta)| \: d\zeta d\bar{\eta} < \infty 
\end{equation}
 We give explicit formulae for the conformal energy of circle homeomorphisms directly in terms of their data. As an example, if $g_o: \mathbb{S} \rightarrow \mathbb{S}$ is an $\eta$-quasi-M\"obius self homeomorphism of the unit circle, then 
\begin{align*}
\mathcal E(g_o) \leq \frac{1}{\pi} \int^{\pi/2}_0 \log \eta \big[ \cot^2(t/2)\big] \,\cos(t) \; dt 
\end{align*}
This estimate is sharp.  Additionally we show how a circle homeomorphism of finite conformal energy can be uniformly approximated on $\IS$ by mappings of strictly smaller energy.
\end{abstract}

\section{Introduction}
The notion of conformal energy was introduced in \cite{AIMO} to study extremal problems for mappings of finite distortion akin to the classical problem of identifying extremal quasiconformal mappings in Teichm\"uller theory.  There (and in sequels \cite{astaladeformations,MM}), interesting new phenomena such as the ``Nitsche phenomena'' were observed in these extremal (or minimisation) problems. There are also connections between conformal energy, which is in a sense dual to the energy introduced by Douglas \cite{douglas1931solution},  and the harmonic mapping theory associated with the Dirichlet extension of self homeomorphisms of the circle.  

Let us begin with a few definitions.  

\subsection{Finite distortion} A mapping $f:\ID\to\ID$ of the unit disk is said to have {\em finite distortion} if 
\begin{enumerate}
\item $f\in W^{1,1}_{loc}(\ID)$,  the space of mappings whose first derivative is locally integrable.
\item $J(z,f)\in L^1_{loc}(\ID)$,  the Jacobian determinant, $J(z,f)={\rm det} \,Df(z)$,  is locally integrable.
\item The distortion
\[ \IK(z,f) :=\frac{\|Df(z)\|^2}{J(z,f) } \]
is finite almost everywhere.
\end{enumerate}
We identify \cite{AIM,IM1} as comprehensive references to the theory of mappings of finite distortion. A mapping is {\em quasiconformal} if it is a homeomorphism and $\IK(z,f)\in L^{\infty}(\ID)$ -- quasiconformal mappings have bounded distortion. 

\subsection{Conformal Energy} 
Let $\IS=\partial\ID$ denote the unit circle and let
$$ g_o:\mathbb{S}\rightarrow\mathbb{S} $$
be a homeomorphism.  We say $g_o$ has finite \textit{conformal energy} if the real number 
\begin{equation}\label{energy} \mathcal E(g_o)=-\frac{1}{2\pi^2}\iint_{\mathbb{S}\times \mathbb{S}}\log |g_o(\zeta)-g_o(\eta)| \: d\zeta d\bar{\eta} < \infty 
\end{equation}
When finite, this integral converges absolutely for homeomorphisms, \cite{AIMO}.

\begin{lemma} Conformal energy is conformally invariant:  If $\phi:\bar\ID\to\bar\ID$ is M\"obius,  then 
\[ {\mathcal E}(\phi|\IS \circ g_o) = {\mathcal E}(g_o)\]
\end{lemma}
\noindent{\bf Proof.} We have 
\[\phi(w)=\zeta \frac{w-a}{1-\bar a w},\quad |a|<1, \; |\zeta|=1 \]
so that we must establish
\[ \iint_{\mathbb{S}\times \mathbb{S}}\log \Big|\frac{g_o(\zeta)-a}{1-\bar a g_o(\zeta)}-\frac{g_o(\eta)-a}{1-\bar a g_o(\eta)}  \Big| \: d\zeta d\bar{\eta} = \iint_{\mathbb{S}\times \mathbb{S}}\log |g_o(\zeta)-g_o(\eta)| \: d\zeta d\bar{\eta} \]
We observe that
\[  \frac{g_o(\zeta)-a}{1-\bar a g_o(\zeta)}-\frac{g_o(\eta)-a}{1-\bar a g_o(\eta)}  = \frac{g_o(\zeta)-g_o(\eta)}{(1-\bar a g_o(\zeta))(1-\bar a g_o(\eta))},  \] 
and so the result follows once we observe Fatou's theorem implies
\begin{eqnarray*}
\lefteqn{ \iint_{\mathbb{S}\times \mathbb{S}} \log |1-\bar a g_o(\eta)| \: d\zeta d\bar{\eta} } \\
&=&\int_{\mathbb{S}} d\zeta \times \int_{\mathbb{S}} \log |1-\bar a g_o(\eta)| \: d\bar{\eta} =0
\end{eqnarray*} 
as the first term vanishes and the second is bounded as $1-|a|\leq |1-\bar a g_o(\eta)|\leq 1+|a|$.  This similarly holds for the remaining term.  
 \hfill $\Box$

  \medskip
  The following corollary is straightforward by applying the above lemma to the identity and making a calculation (which we do not do as we essentially reproduce it in a more general setting later).
  \begin{corollary}\label{corollary1}  If $\phi:\IS\to\IS$ is M\"obius,  then 
$ {\mathcal E}(\phi) = 1$ \end{corollary}
  The converse to Corollary \ref{corollary1} is true as we shall see below.  But it does not seem to be an easy consequence of the definition.
  
  \medskip
  
  In a similar vein we note the following.  
We say that $f:\IS\to \IS$ is bilipschitz if there is a constant $L \geq 1$ such that 
$$\frac{1}{L}|\zeta - \eta| \leq |f(\zeta)-f(\eta)| \leq L|\zeta - \eta|$$
Then we have the following Lipschitz invariance of conformal energy.

\begin{lemma}  If $f:\IS\to\IS$ is $L$-bilipschitz, and $g_o:\IS\to\IS$,  then 
\[\max\big\{1,{\mathcal E}(g_o)-\frac{1}{\pi^2} \log L\big\} \leq {\mathcal E}(f \circ g_o) \leq {\mathcal E}(g_o)+\frac{1}{2\pi^2} \log L\]
and consequently 
\[  {\mathcal E}(f ) \leq  \frac{1}{2\pi^2} \log L \]
\end{lemma}

Although an elementary calculation,   the logarithmic term is a little surprising as we will see later.

\section{Extremal problems}

We now recall the main theorem of \cite[Theorem 11.4]{AIMO}:
\begin{theorem*} \label{thm1}
Let $g_o: \mathbb{S} \rightarrow \mathbb{S}$ be a homeomorphism.  Consider the minimisation problem to identify
\begin{equation}  \inf_{g} \left\{ \frac{1}{\pi} \iint_{\mathbb{D}} \mathbb{K}(z,g)\; |dz|^2 \right\}
\end{equation}
where the infimum is taken over all homeomorphisms $g:\overline{\mathbb{D}} \rightarrow \overline{\mathbb{D}}$ of finite distortion for which $g|\mathbb{S}=g_o$.  Then there exists a  minimiser if and only if $g_o$ has finite energy. This minimum value is $\mathcal E(g_o)$,  and the minimiser is unique.
\end{theorem*}

\medskip

Note therefore,  as a particular consequence, if $g_o$ has a quasiconformal extension to the disk, with say $\IK=\|\IK(z,f)\|_{L^\infty(\ID)}$, then  
$$1 \leq \mathcal E(g_o) \leq \mathbb{K}$$

We now have the following converse to Corollary \ref{corollary1}.

\begin{corollary} \label{cor1}
$\mathcal E(g_o)=1$ only if $g_o$ represents the boundary values of a M\"obius transformation.
\end{corollary}
\noindent{\bf Proof.} First observe that $\IK(z,f)\geq 1$.  Then $E(g_o)=1$ implies $\IK(z,f)=1$ almost everywhere.  As $f$ is a homeomorphism and $J(z,f)=|f_z|^2-|f_\zbar|\geq 0$ almost everywhere, we have
\[ \IK(z,f) = \frac{|f_z|^2+|f_\zbar|^2}{|f_z|^2-|f_\zbar|^2}=1\Rightarrow f_\zbar \equiv 0 \;\; \mbox{ almost every $z\in\ID$}.\]
The equation $f_\zbar=0$ with $f\in W^{1,1}_{loc}(\ID)$ implies $f$ is conformal by the Looman-Menchoff theorem, see \cite{Nara}. \hfill $\Box$

\section{Quasi-M\"obius mappings}

Quasi-M\"obius mappings of fairly general  spaces were introduced by J. V\"ais\"al\"a \cite{V} extending his earlier work. with P. Tukia \cite{TV}.  Since then these mappings have found wide application in geometric analysis,  particularly on metric spaces.

\medskip

A homeomorphism of the circle is said to be $\eta$-quasi-M\"obius (where $\eta:[0,  \infty]\rightarrow [0, \infty]$ is increasing with $\eta(0)=0$ and $\eta(\infty)=\infty$) if for every quadruple of points 
$$\frac{|g_o(a)-g_o(b)||g_o(c)-g_o(d)|}{|g_o(a)-g_o(c)||g_o(b)-g_o(d)|} \leq \eta \left( \frac{|a-b||c-d|}{|a-c||b-d|} \right)$$
Interchanging the roles of $a$ and $b$ 
immediately implies the two sided inequality
$$\frac{1}{\eta(1/[a,b,c,d])} \leq [g_o(a),g_o(b),g_o(c),g_o(d)] \leq \eta([a,b,c,d])$$
where we have written the cross ratio:
$$[a,b,c,d]=\frac{|a-b||c-d|}{|a-c||b-d|}.$$
Thus $1\leq \eta(t)\eta(1/t)$ and $\eta(1) \geq 1$.  

\medskip

We recall the following known facts from \cite{V}. A bi-Lipschitz map is quasisymmetric and quasi-M\"obius. A quasisymmetric map is quasi-M\"obius,  and finally a quasi-M\"obius map defined on a bounded space is quasisymmetric. In particular,  it follows that if $g:\ID\to\ID$ is $K$-quasiconformal,  then $f$ has a unique boundary value  extension $f_o:\IS\to\IS$ which is $\eta$-quasi-M\"obius.  Here $\eta$ depends on $K$ within bounds,  but the exact relationship is not known.  The next result is one of our main theorems.

\begin{theorem}\label{main}
Let $g_o: \mathbb{S} \rightarrow \mathbb{S}$ be $\eta$-quasi-M\"obius.  Then 
\begin{align*}
 \mathcal E(g_o) \leq \frac{1}{\pi} \int^{\pi/2}_0  \log \big[ \eta (\cot^2(t/2))\big]\,\cos(t)   \,  dt
\end{align*}
This estimate is sharp.
\end{theorem}
\noindent{\bf Proof.}   We first observe
\begin{align*}
\mathcal E(g_o) &= -\frac{1}{2\pi^2}\iint_{\mathbb{S}\times\mathbb{S}} \log |g_o(\zeta)-g_o(\xi)|\: d\zeta d\bar{\xi}  
 =\frac{1}{2\pi^2}\iint_{\mathbb{S}\times\mathbb{S}} \log |g_o(\zeta)-g_o(-\xi)|\: d\zeta d\bar{\xi} \\
&=-\frac{1}{2\pi^2}\iint_{\mathbb{S}\times\mathbb{S}} \log \frac{1}{|g_o(\zeta)-g_o(-\xi)|}\:d\zeta d\bar{\xi} 
\end{align*}
Therefore
\begin{align*}
4&\mathcal E(g_o) = \\
&=\frac{1}{2\pi^2}\iint_{\mathbb{S}\times\mathbb{S}} \log \frac{|g_o(\zeta)-g_o(-\xi)||g_o(-\zeta)-g_o(\xi)|}{|g_o(\zeta)-g_o(\xi)||g_o(-\zeta)-g_o(-\xi)|}\: d\zeta d\bar{\xi}\\
&=\frac{1}{2\pi^2}\int^{2\pi}_0\int^{2\pi}_0 \log g_o[e^{it}, e^{is}, -e^{it},-e^{is}]\cos(t-s)\:dtds \\
&=\frac{1}{2\pi^2}\int^{2\pi}_0\int^{2\pi}_0 \log g_o[e^{i(t+s)}, e^{is}, -e^{i(t+s)},-e^{is}]\cos(t)\:dtds
\end{align*}
Next we observe that
$$[e^{i(t+s)}, e^{is}, -e^{i(t+s)},-e^{is}]=[e^{it}, 1, -e^{it},-1]=\cot^2(t/2).$$

We independently considering where $\cos(t)$ is positive and negative the energy bound is proved.  Sharpness is established once we observe that if $\eta(t)=t$,  so that $g_o$ are the boundary values of a M\"obius transformation, then 
\begin{equation}\label{calc} {\mathcal E}(g_o)=\frac{1}{\pi }\int_0^{\frac{\pi }{2}} \cos (t) \log \big[\cot ^2(t/2)\big] \, dt = 1 
\end{equation}
\hfill $\Box$

\medskip
We give further examples below at \S\ref{examples} to show that the estimate of Theorem \ref{main} is of the correct order.

\medskip
\subsection{$\eta(t)=\alpha t$, $t\geq 1$} Here we note an elementary consequence of the last calculation above for mappings that distort cross ratios by a bounded amount.
\begin{theorem} Suppose $g_o:\IS\to\IS$ is a homeomorphism and that there is $\alpha\geq 1$ so that $g_o$ has bounded cross-ratio distortion:
$$\frac{1}{\alpha}\, [a,b,c,d] \leq [g_o(a),g_o(b),g_o(c),g_o(d)] \leq \alpha\, [a,b,c,d]. $$ 
Then 
\[ {\mathcal E}(g_o)\leq 1+\frac{1}{\pi} \; \log(\alpha).\]
\end{theorem}
\noindent{\bf Proof.} This follows from Theorem \ref{main} with $\eta(t)=\alpha \,t$ and the calculation at (\ref{calc}).\hfill $\Box$

\section{Some examples}\label{examples} In this section we give a couple of relevant examples.

First, define a homeomorphism $\theta:[0,2\pi]\rightarrow [0,2\pi]$ by
\begin{equation}
 \theta(t) = 
\begin{cases}
t/ \lambda,  & 0\leq t \lambda \\
1+(2\pi-1)/(2\pi -\lambda)t, & \lambda \leq t \leq 2\pi,
\end{cases} 
\end{equation}
and then define a homeomorphism of the circle as 
\begin{equation}
f_o(e^{it})=e^{i\theta(t)}.
\end{equation}
\begin{theorem} The energy $\mathcal E(f_o)$ is bounded independently of $\lambda$ while the energy of the inverse $h_o=f_o^{-1}:\IS\to\IS$ has  $\mathcal E(h_o)$.
\end{theorem}
We compute the energy of the map $f_o$. Set 
\[ \alpha =(2\pi-1)/(2\pi -\lambda), \quad F(t,s)=\cos(t-s)\, \log |e^{i\theta(t)}-e^{i\theta(s)}|.\]

With $\zeta=e^{it}, \eta=e^{is}$ we see (using symmetry),
\begin{align*}
&\iint_{\mathbb{S}\times\mathbb{S}} \log |f_o(\zeta)-f_0(\eta)|\:d\zeta d\bar{\eta} =\int^{2\pi}_0\int^{2\pi}_0 \log|e^{i\theta(t)}-e^{i\theta(s)}|\cos(t-s)\:dtds \\
&=\int^{\lambda}_0\int^{\lambda}_0F(s,t)dtds+\int^{2\pi}_{\lambda}\int^{2\pi}_{\lambda}F(s,t)dtds+2\int^{\lambda}_0\int^{2\pi}_{\lambda}F(t,s)\:dtds
\end{align*}
We will now deal with each of these three integrals in turn.
\begin{align*}
&\int^{\lambda}_0 \int^{\lambda}_0F(s,t)\:dtds = \int^{\lambda}_0\int^{\lambda}_0 \log |1-e^{i(t-s)/\lambda}|\cos(t-s)\:dtds \\
&= \lambda\int^{\lambda}_0\int^{1-s/\lambda}_{-s/\lambda} \log |1-e^{iu}|\cos(\lambda u)\:duds \leq \lambda^2 \int^{2\pi}_0 \left| \log |1-e^{iu}| \right| du
\end{align*}
where we have used the periodicity of the integrand.  As $\log |1-e^{iu}| \in L^1[0,2\pi]$ this integral tends to 0 with $\lambda$:
\begin{equation}
\int^{\lambda}_0\int^{\lambda}_0 F(s,t)\:dtds\rightarrow 0 \: \text{as} \: \lambda \rightarrow 0
\end{equation}
Next,
\begin{align*}
&\int^{2\pi}_{\lambda}\int^{2\pi}_{\lambda}F(s,t)\:dtds = \int^{2\pi}_{\lambda}\int^{2\pi}_{\lambda} \log |1-e^{i\alpha(t-s)}|\cos(t-s)\:dtds \\
&= \int^{2\pi}_{\lambda}\int^{\alpha(2\pi-s)}_{\alpha(\lambda-s)} \log|1-e^{iu}|\cos(u/\alpha)\:dtds
\end{align*}
Here we can apply the dominated convergence theorem to see 
\begin{equation}
\int^{2\pi}_{\lambda}\int^{2\pi}_{\lambda}F(s,t)\:dtds \rightarrow \int^{2\pi}_0 \int^{a(2\pi-s)}_{-as} \log |1-e^{i(t-s)}|\cos(a(t-s))\:dtds
\end{equation}
as $\lambda \rightarrow 0$ and $a=2\pi /(2\pi -1)$.  This last number can be computed numerically.  
Finally we have to deal with the cross terms,  namely
\begin{align*}
&\left| \int^{\lambda}_0\int^{2\pi}_{\lambda} F(t,s)\:dtds \right| =\left| \int^{\lambda}_0\int^{2\pi}_{\lambda}\log|1-e^{i(1+\alpha t-s/\lambda)}|\cos(t-s)\:dtds \right| \\
&\leq \frac{\lambda}{\alpha}\int^1_0\int^{1+2\pi\alpha}_{1+\alpha\lambda}\left| \log|1-e^{i(u-v)}| \right|\: dudv
\end{align*}
As above we deduce 
\begin{equation}
\int^{\lambda}_0 \int^{2\pi}_{\lambda}F(t,s)\:dtds \rightarrow 0 \: \text{as} \: \lambda \rightarrow 0
\end{equation}
When we deal with the inverse function the analysis is similar.  However is it clear that the behaviour of the integral is dominated by the term
\begin{align*}
&\int^1_0 \int^1_0 \log |e^{i\lambda t}-e^{i\lambda s}| \cos(t-s) \:dtds =\int^1_0 \int^1_0 \log|1-e^{i\lambda (t-s)}|\cos(t-s)\:dtds \\
&=\int^1_0 \int^1_0 \log|\lambda(t-s)[1+h(\lambda(t-s))]|\cos(t-s) \:dtds
\end{align*}
Here $h$ is the remainder term in the Taylor series expansion and can be assumed uniformly small if we choose $\lambda$ small enough.  Thus
\begin{align*}
&\int^1_0\int^1_0 \log |e^{i\lambda t}-e^{i\lambda s}| \cos(t-s)\:dtds \\
&= |\log \lambda| \int^1_0 \int^1_0 \cos(t-s)dtds+\int^1_0\int^1_0\log |(t-s)|\cos(t-s)\:dtds +O(1) \\
&\approx (2-2\cos(1))|\log \lambda| +O(1)
\end{align*}
which diverges with $\lambda$.  The constants associated with the definition of quasi-M\"obius are estimated by considering the four points
$$[1,e^{i\lambda}, e^{i\pi},e^{3i\pi/2}]=\frac{|1-e^{i\lambda}||e^{i\pi}-e^{3i\pi/2}|}{|1-e^{i\pi}||e^{i\lambda}-e^{3i\pi/2}|}\approx \frac{\lambda}{2}$$
The image cross ratio is 
$$[1,e^i, e^{ i\alpha\pi}, e^{3i\alpha \pi/2}]=\frac{|1-e^i||e^{i\alpha\pi} - e^{3i\alpha\pi/2}|}{|1-e^{i\alpha \pi}||e^i-e^{3i\alpha\pi/2}|} \approx 0.541 \ldots$$
So
$$f_o[1,e^{i\lambda}, e^{i\pi},e^{3i\pi/2}]=C\, \lambda\, [1,e^{i\lambda}, e^{i\pi},e^{3i\pi/2}]$$
with $C\approx 1.08 \ldots$.  In particular, for any function $\eta$ for which $f_o$ is $\eta$-quasi-M\"obius we have
$$\eta(2) \geq \frac{1}{\lambda}$$
Clearly the bilipshitz constant is $1/\lambda$.  The map $f_o$ has finite energy independent of $\lambda$ while its inverse has energy $\approx \log \lambda$.  This is aligned with our more general calculations earlier and shows those estimates to be of the correct orders. 

\subsection*{A Map and Inverse of Finite Energy but not Quasim\"obius}
A quasi-M\"obius map always has finite energy and since the inverse of a quasi-M\"obius map is again quasi-M\"obius, the inverse will also have finite energy.  We next give an example to show this does not characterise quasi-M\"obius maps.  

\medskip

Consider the map $\theta : [0, 2\pi] \rightarrow [0,2\pi]$ by
\begin{equation}
\theta(t) =
\begin{cases}
t^2, & 0\leq t \leq 1 \\
t, & 1\leq t \leq 2\pi
\end{cases}
\end{equation}
and define a homeomorphism of the circle as
\begin{equation}
g_o(e^{it})=e^{i\theta(t)}
\end{equation}
We leave it to the reader to verify the elementary fact that the the mapping $g_o$ is not quasi-M\"obius.

Next,  following the calculations above,  we will see that $g_o$ and $g_o^{-1}$ have finite energy.We first examine the integrals
$$\int^1_0 \int^1_0 \left| \log |e^{it^2}-e^{is^2}| \right| dtds \quad \text{and} \quad \int^1_0 \int^1_0 \left| \log |e^{i\sqrt{t}}-e^{i\sqrt{s}}| \right| dtds$$
As above we consider the power series expansion of the exponential function to see 
\begin{align*}
\int^1_0  &\int^1_0 \left| \log |e^{it^2} -e^{is^2}| \right|dtds = \int^1_0 \int^1_0 \left| \log |t^2-s^2| \right| dtds +O(1) \\
& \leq \int^1_0 \int^1_0 \left| \log |t-s| \right| dtds +\int^1_0 \int^1_0 \left| \log |t+s| \right| dtds + O(1) \\
& \leq 2\int^1_0 \int^1_0 \left| \log |t-s| \right| dtds +O(1)
\end{align*}
and this last integral is finite.  Further
\begin{align*}
\int^1_0 & \int^1_0 \left| \log |e^{i\sqrt{t}} -e^{i\sqrt{s}}| \right|dtds \leq \int^1_0 \int^1_0 \left| \log |\sqrt{t} -\sqrt{s}| \right|dtds +O(1) \\
&\leq \int^1_0 \int^1_0 \left| \log |t -s| \right|dtds +O(1)
\end{align*}
is again finite.  A typical cross term
$$\int^1_0 \int^{2\pi}_1 \left| \log |e^{it^2}-e^{is}| \right| dtds = \int^1_0 \int^{2\pi}_1 \left| \log |t^2-s| \right| dtds +O(1)$$
and $|t-s| \leq |t^2-s| \leq 2\pi$.  So this integral is dominated by those above.  It now follows both $g_o$ and $g_o^{-1}$ have finite energy. \hfill $\Box$

\section{Critical Points}
In this section we consider the critical points of the energy functional.  A homeomorphism $g_o:\mathbb{S} \rightarrow \mathbb{S}$ can be written in the form 
$$g_o(e^{i x})=e^{i\theta(x)},  \: x \in [0,  2\pi]$$
After some normalisations,  $\theta : [0,2\pi]\rightarrow [0,2\pi]$ may be assumed to be a strictly increasing continuous function $\theta(0) = 0$,  $\theta(2\pi)=2\pi$.  A variation of $g_o$ can be obtained by considering a smooth function $\phi :[0,2\pi]\rightarrow \mathbb{R}$ with $\phi(0)=\phi(2\pi)=0$, $|\nabla \phi|\leq 1$. We set
$$g_o^t(e^{ix}) = e^{i(\theta(x)+ t \phi(x))}, \: x \in [0,2\pi],  \quad g_o^t\Big|_{t=0}=g_o$$
Note that $g_o^t:\IS\to\IS$ need not be a homeomorphism of the circle,  though it is a map $\mathbb{S}\rightarrow \mathbb{S}$ and one can still speak of its conformal energy.  In fact if $\theta'\geq\epsilon>0$ is bounded below, then for all small $t$ we will obtain a homeomorphism as $\theta(x)+i t \phi(x)$ will be increasing and $2\pi$-periodic.  For simplicity we will write $\mathcal E(\theta + t\phi)$ for $\mathcal E(e^{i(\theta + t\phi)})$.  Further,  we extend $\theta$ and $\phi$ periodically to the real line and consider
\begin{align*}
&\frac{d}{dt} \mathcal E(\theta +t\phi) \\
&= \frac{1}{2}\, \frac{d}{dt} \left\{ \int^{2\pi}_0 \int^{2\pi}_0 \log \left|e^{i(\theta(x) + t\phi(x))}-e^{i(\theta(y)+t\phi(y))} \right|^2 \cos(x-y) \:dxdy \right\} \end{align*}
Then computing this derivative and evaluating at $t=0$ gives
\begin{align} 
&\frac{d}{dt} \mathcal E (\theta +t\phi)\Big|_{t=0} \nonumber \\ 
&= \frac{1}{2}\int^{2\pi}_0 \int^{2\pi}_0 \cot \left[ \frac{\theta(x)-\theta(y)}{2} \right] (\phi(x)-\phi(y)) \cos(x-y) dxdy \label{t=0}
\end{align}
Here we must discuss the convergence of this integral.  The function $\phi$ is smooth and so near the diagonal $\{ x=y \}$ where the singularity occurs we have $$\phi(x)-\phi(y) = \phi '(s)(x-y)$$ for some $s \in [x,y]$.  As $\phi'$ is bounded, we really must assume that
\begin{equation} \label{eq12}
\int^{2\pi}_0 \int^{2\pi}_0 \left| \frac{x-y}{\theta(x)-\theta(y)} \right| dxdy <\infty
\end{equation}
 To move forward,  we therefore make the assumption that there are constants $\alpha$ and  $p<2$ such that 
\begin{equation}
|\theta(x)-\theta(y)| \geq \alpha|x-y|^p
\end{equation}
Now at least we can discuss convergent integrals. Using $2\pi$-periodicity and a linear change of variables $(x\rightarrow x+y)$, 
\begin{align*}
&\frac{d}{dt} \mathcal E(\theta+t\phi) \Bigg|_{t=0} \\
&=\frac{1}{2}\int^{2\pi}_0 \int^{2\pi}_0  \frac{\sin(\theta(x+y)-\theta(y))}{1-\cos(\theta(x+y)-\theta(y))}(\phi(x+y)-\phi(y))\cos(x)\:dxdy
\end{align*}
Now for each $x \in [0,2\pi]$ the integral
$$F(x)=\int^{2\pi}_0 \frac{\sin(\theta(x+y)-\theta(y))}{1-\cos(\theta(x+y)-\theta(y))}(\phi(x+y)-\phi(y))dy$$
converges and with our assumptions defines a continuous bounded function of $x$.  Similarly, 
$$G(y)=\int^{2\pi}_0 \frac{\sin(\theta(x+y)-\theta(y))}{1-\cos(\theta(x+y)-\theta(y))}(\phi(x+y)-\phi(y))\cos(x)dx$$
defines a continuous and bounded function of $y$.  Set
$$A(x)=\int^{2\pi}_0 \frac{\sin(\theta(x+y)-\theta(y))}{1-\cos(\theta(x+y)-\theta(y))}\phi(x+y)dy$$
$$B(x)=-\int^{2\pi}_0 \frac{\sin(\theta(x+y)-\theta(y))}{1-\cos(\theta(x+y)-\theta(y))}\phi(y)dy$$
so that for $x \in (0, 2\pi)$ we have $F(x) = A(x)+B(x)$.  Of course as $x \rightarrow 0$, both $A(x)$ and $B(x)$ tend to $\infty$ (mod $2\pi$).  Using periodicity again we have 
$$B(x)=-\int^{2\pi}_0 \frac{\sin(\theta(2x+y)-\theta(x+y))}{1-\cos(\theta(2x+y)-\theta(x+y))}\phi(x+y)dy$$
Putting these back together gives
\begin{align*}
&F(x) \\
&= \int^{2\pi}_0 \left( \frac{\sin(\theta(x+y)-\theta(y))}{-\cos(\theta(x+y)-\theta(y))} - \frac{\sin(\theta(2x+y)-\theta(x+y))}{1-\cos(\theta(2x+y)-\theta(x+y))} \right) \phi(x+y) dy \\
&= \int^{2\pi}_0 \left( \frac{\sin(\theta(u)-\theta(u-x))}{-\cos(\theta(u)-\theta(u-x))} - \frac{\sin(\theta(u+x)-\theta(u)}{1-\cos(\theta(u+x)-\theta(u))} \right) \phi(u) du
\end{align*}
Then
\begin{align*}
&\int^{2\pi}_0F(x)\cos(x)dx \\
&= \int^{2\pi}_0 \left[ \int^{2\pi}_0 \left( \cot \frac{\theta(u)-\theta(u-x)}{2} - \cot \frac{\theta(u+x)-\theta(u)}{2} \right)\phi(u)du \right] \cos(x) dx
\end{align*}
Fubini's theorem allows us to change the order of integration here provided
\begin{equation} \label{eq114}
\cot\frac{\theta(u)-\theta(u-x)}{2} - \cot\frac{\theta(u+x)-\theta(u)}{2}
\end{equation}
is integrable.  Let us assume this for now and discuss the conditions later.  This leads us to the formula
\begin{align*}
&\frac{d}{dt}\mathcal E(\theta+t\phi)\Bigg|_{t=0} \\
&= \frac{1}{2} \int^{2\pi}_0 \int^{2\pi}_0 \left( \cot\Big[\frac{\theta(y)-\theta(y-x)}{2}\Big] - \cot\Big[\frac{\theta(y+x)-\theta(y)}{2} \Big]\right) \cos(x)dx \phi(y)dy
\end{align*}
It further follows from Fubini's theorem together with the fact that $\phi$ is an essentially arbitrary function and so we see that at a critical point of the conformal energy (homeomorphism $e^{i\theta(x)}$) we must have 
$$\frac{d}{dt} \mathcal E(\theta + t\phi) \Bigg|_{t=0} =0$$
and thus for every $y \in [0,2\pi]$ the integral
\begin{equation}\label{16}
\int^{\pi}_{-\pi} \left( \cot \Big[\frac{\theta(y)-\theta(y-x)}{2} \Big]-\cot\Big[\frac{\theta(y+x)-\theta(y)}{2}\Big] \right) \cos(x)dx=0
\end{equation}
The convergence of this integral will negate our postponed discussion of the convergence of the integral at \eqref{eq114}.  If $\theta$ is $C^{1, \alpha}(\mathbb{S})$ for any $\alpha >0$ then this integral is finite. 

\begin{theorem} A sufficiently regular critical point,  for instance $C^{1, \alpha}(\mathbb{S})$ for any $\alpha >0$, of the energy satisfies (\ref{16})
for all $x\in[-\pi,\pi]$.
\end{theorem}
\medskip

If we write 
\begin{equation}
\theta(t)=2\arctan[u(t)],  \quad \theta(t)=\arctan[v(t)],\;\;  v=\frac{2u}{1-u^2}.
\end{equation}
We may do that for an increasing $u$. Then the integrand here becomes
\[ \frac{u(y)^2+1}{u(y)-u(y-x)}+\frac{u(y)^2+1}{u(y)-u(x+y)}-2 u(y) \]
Thus (\ref{16}) is equivalent to
\begin{equation}
\int^{\pi}_{-\pi} \left( \frac{u(y)^2+1}{u(y)-u(y-x)}+\frac{u(y)^2+1}{u(y)-u(x+y)}-2 u(y) \right) \cos(x)dx=0
\end{equation}
and hence
\begin{equation}
\int^{\pi}_{-\pi} \left(\frac{1}{u(x+y)-u(y)} - \frac{1}{u(y)-u(y-x)} \right) \cos(x)dx=0
\end{equation}
This seems closely related to the (reciprocal of the) symmetric derivative of $u$.
 
\begin{conj} \label{conj1}
If $\theta$ is a sufficiently regular homeomorphism $[0,2\pi] \rightarrow [0,2\pi]$ with $\theta(\pi)=\pi$,  and if for every $y \in [0,2\pi]$ the integral
\begin{equation}
\int^{2\pi}_0 \left( \cot\Big[ \frac{\theta(y)-\theta(y-x)}{2}\Big] -\cot\Big[\frac{\theta(y+x)-\theta(y)}{2}\Big] \right) \cos(x)dx=0
\end{equation}
then there is $a>0$ such that 
$$\theta(t) =\arctan \left( \frac{(1-a^2)\sin(t)}{(1+a^2)\cos(t)-2a} \right)$$
\end{conj}
We note that if $f(z)=(z-a)/(1-az)=e^{i\theta(t)}$,  $a \in \mathbb{R}$  is a M\"obius transformation fixing $\pm 1$, the formula for $\theta$ is 
\begin{align*}
\theta(t) &= \arg ((e^{it}-a)/(1-ae^{it})) \\
&= \arg (e^{it}-2a+a^2e^{-it}) \\
&= \arg ((1+a^2) \cos(t) -2a + i(1-a^2)\sin(t)) \\
&= \arctan \left( \frac{(1-a^2)\sin(t)}{(1+a^2)\cos(t)-2a} \right)
\end{align*}

\begin{conj} \label{conj2}
M\"obius transformations are the only homeomorphic critical points of the conformal energy functional $\mathcal E: Hom(\mathbb{S}) \rightarrow [1, \infty)$.
\end{conj}

The subtlety here is with the setup we have arranged to get equations -- it is via smooth perturbations and so the conjecture basically asks for a smooth perturbation which lowers energy. If one is not concerned with smooth perturbations then we give the following.

\begin{theorem} Let $g_o:\IS\to\IS$ be a homeomorphism of finite energy which does not represent the boundary values of a M\"obius transformation. Then there is a family of homeomorphisms $\{g^{t}:\IS\to\IS:t\in [0,1)\}$ such that
\begin{enumerate}
\item $g^{t} \to g_o$ uniformly on $\IS$ as $t\to 1$.
\item $\mathcal{E}(g^t) < \mathcal{E}(g_o)$ and $\mathcal{E}(g^t) \to \mathcal{E}(g_o)$ as $t\to 1$.
\item $g^{t}$ converge uniformly to the bounadry values of a M\"obius transformation as $t\to 0$.
\end{enumerate}
\end{theorem}
\noindent{\bf Proof.}   We let $h_o=g_o^{-1}:\IS\to\IS$ and let $H:\ID\to\ID$ be the harmonic extension of $h_o$ to the unit disk.  Then $H$ is a diffeomorphism,  and we set $G=H^{-1}:\ID\to\ID$ \cite{Duren}.  Of course $G$ extends homeomorphically to the boundary $\partial\ID=\IS$ and $G|\IS=g_o$.  As $H$ is harmonic, $H_w=\phi$ and $H_\wbar = \bar\psi$ and $\phi$ and $\psi$ are holomorphic. Then we define (in the notation of P. Duren \cite{Duren}) the $2^{nd}-$Beltrami equation of $H$ as
\begin{equation}\label{21}
H_\wbar = \overline{\nu(w)}\; \overline{H_w}.
\end{equation}
As $H$ is a diffeomorphism, $H_w=\phi\neq0$ and so $\nu$ is holomorphic as well as  $|\nu(w)|< 1$, $w\in \ID$.
From (\ref{21}) we see that $G$ solves the equation 
\begin{equation}\label{22}
G_\zbar = -\overline{\nu(G)}\; G_z.
\end{equation}
Thus the Beltrami coefficient for $G$ satisfies  
\[\mu(z)= G_\zbar/G_z, \quad |\mu(z)|=|\nu(G)|<1,   z\in \ID \]
Also $G_z$ does not vanish,  so $\mu$ is smooth in $\ID$.  Also,  by hypothesis $\mu$ is not identically $0$.
 For $t\in [0,1)$ let $G^t:\ID\to\ID$ be the (quasiconformal) solution to the (uniformly elliptic) Beltrami equation
 \begin{equation}
 G^t_\zbar = t \mu(z) G^t_z, \;\; G^t(0)=G(0),  \quad z\in\ID.
 \end{equation}
and set $g^t = G^t|\IS$. Now \cite{AIMO} recognises $G$  as being extremal for its boundary values and so we now make the following calculation using Theorem \ref{thm1}
\begin{eqnarray*}
\mathcal{E}(g^t) &\leq & \frac{1}{\pi} \int_\ID \IK(z,G^t) \; dz  = \frac{1}{\pi} \int_\ID \frac{1+t|\mu|^2}{1-t|\mu|^2}\, dz 
\end{eqnarray*}
Next,
\[ \frac{d}{dt} \frac{1}{\pi} \int_\ID \frac{1+t|\mu|^2}{1-t|\mu|^2}\, dz = \frac{1}{\pi} \int_\ID \frac{2|\mu|^2}{1-t|\mu|^2}\, dz \geq \frac{1}{\pi} \int_\ID \frac{2|\mu|^2}{1-|\mu|^2}\, dz  > 0,\]
as $\mu$ is not identically $0$.  This says that the energy bound is strictly increasing,  so 
\[ \mathcal{E}(g^t) < \mathcal{E}(g_o), \quad t\in[0,1) \]
and $H^t=(G^t)^{-1}$ have a uniform bound in $W^{1,2}(\ID)$ as 
\[ \int_\ID \|DH^t\|^2 \, dw = \int_\ID \IK(z,G^t) \, dz \leq \int_\ID \IK(z,G) \, dz <\infty \]
Hence,  as homeomorphisms, there is a uniform modulus of continuity and recall that $G^t(0)=G(0)$ for all $t$.  Thus there exists $H^t\to F:\ID\to\ID$, a homeomorphism.  However the family $\{G^t\}$ is locally uniformly quasiconformal and so $G^t\to G$ locally uniformly in $\ID$ and it follows that $F=G^{-1}$,  and hence $G^t\to G$ uniformly in $\ID$.  Thus
\begin{equation}
g^t \to g_o, \quad\mbox{uniformly on $\IS$}.
\end{equation}
Further, we also have $\mathcal{E}(g^t) \to \mathcal{E}(g_o)$ and $\mathcal{E}(g^t) \leq \mathcal{E}(g_o)$. And finally, as $t\to 0$, $g^t$ converges to the boundary values of a M\"obius transformation. \hfill $\Box$

\end{document}